\documentclass[onecolumn]{article}

\usepackage{amsfonts,amsmath,amsthm}

\newtheorem{prop}{Proposition}

\newcommand{\RR}{{\mathbb R}}

\newcommand{\CC}{{\mathbb C}}

\newcommand{\dv}[2]{{\frac{\partial #1}{\partial #2}}}

\newcommand{\dotex}{{\frac{d}{dt}}}

\newcommand{\LL}{{\mathcal L}}
\newcommand{\HH}{{\mathcal H }}
\newcommand{\LLc}{\widetilde{\mathcal L}}
\newcommand{\HHc}{\widetilde{\mathcal H}}

\newcommand{\LLIM}{{\mathcal L}_{\text{\tiny IM}}}
\newcommand{\is}{\imath_{s}}
\newcommand{\ib}{\bar{\imath}}
\newcommand{\ir}{\imath_r}

\begin{document}

\title{Euler-Lagrange models with complex currents of    three-phase electrical machines and observability issues}

\author{
D.~Basic \and
F.~Malrait \and
        P.~Rouchon
\thanks{Duro Basic   and Fran\c{c}ois Malrait are with Schneider Electric, STIE, 33, rue Andr\'{e}
Blanchet, 27120 Pacy-sur Eure. E-mail: duro.basic@schneider-electric.com,
francois.malrait@schneider-electric.com}
\thanks{Pierre  Rouchon is with  Mines ParisTech, Centre Automatique et Syst\`emes, Math\'{e}matiques et Syst\`{e}mes,  60 Bd Saint-Michel, 75272 Paris cedex 06, France.
E-mail: pierre.rouchon@mines-paristech.fr} }

\maketitle

\begin{abstract}

A new  Lagrangian formulation with  complex currents  is developed and yields a direct and simple method for modeling three-phase permanent-magnet and  induction machines. The Lagrangian is the sum a  mechanical one  and of a magnetic one. This magnetic Lagrangian is expressed in terms of rotor angle, complex stator and rotor currents. A complexification procedure widely used in quantum electrodynamic  is applied here in order to derive the Euler-Lagrange equations with  complex  stator and rotor currents. Such complexification process avoids the usual separation into real and imaginary parts and simplifies notably the calculations. Via simple modifications of such  magnetic Lagrangians we derive new  dynamical models describing  permanent-magnet machines with both saturation and saliency, and  induction machines with both magnetic saturation and space harmonics. For each model we also provide its Hamiltonian thus its magnetic energy. This energy is also expressed with complex currents and can be directly used  in Lyapunov and/or passivity based control.  Further, we briefly investigate the observability of this class of  Euler-Lagrange models, in the so-called sensorless case when the measured output is the stator current and the load torque is constant but unknown. For all the dynamical models obtained via such variational principles, we prove that their linear tangent systems are unobservable around  a one-dimensional family of steady-states attached to the same  constant stator voltage and current. This negative result explains why sensorless control of three-phase electrical machines  around zero stator frequency remains yet a difficult control problem.

\end{abstract}


\section{Introduction}

Modeling electrical machines  with  magnetic-saturation and space-harmonics effects   is not a  straightforward task
and could lead to complicated  developments when a  detailed physical  description is  included  (see,
e.g.,\cite{chiasson:book05,InductionMotor-handbook}). Even if such effects are not dominant  they could  play an important
role for sensorless control (no rotor position or velocity  sensor). For a permanent-magnet
machine, the rotor position will be   unobservable without saliency. For standard models of induction machines (no magnetic saturation, no space-harmonics) the rotor velocity  is always unobservable at zero stator frequency~\cite{canudas-et-al:cdc00,holtz-procIEEE:02}.
A global observability analysis based on such standard models is given in~\cite{ibarra-el-al:auto04}. In this note we develop  a systematic method to include into such standard control models  magnetic-saturation, saliency and/or space-harmonics effects. Our initial motivation  was to see whether such  non-observability is destroyed by such modelling changes or not. It appears that any physical-consistent model admits the same kind of observability deficiency at zero stator frequency.

By physically consistent models we mean  Lagrangian-based models. We propose here an extension of  Lagrangian modeling of three-phase machines  with  real variables (see, e.g.\cite{ortega-et-al-book98})  to complex electrical variables. It is directly inspired from quantum electro-dynamics where Lagrangian with complex generalized positions and velocities  are widely used (see, e.g.,~\cite{cohenT-book1}, page 87). We obtain, from such  Lagrangian  functions,   physically consistent and synthetic  Euler-Lagrange  models   directly
expressed  with  complex stator and rotor currents. Such modeling method   by-passes  the usual detailed physical descriptions that are not easily accessible to the control community. Here we propose a much more direct way: it just consists in modifying the magnetic part of the Lagrangian directly  expressed with complex currents and then in  deriving the dynamic equations from the Euler-Lagrange with complex variables. We obtain automatically the dynamics of the electrical part as a set of complex differential  equations. We suggest here  simple Lagrangians modeling simultaneously magnetic-saturation, saliency and space-harmonic effects. The obtained dynamics extend directly the ones used  in almost all control-theoretic papers and include also  more elaborate ones that can be found in specialized  books such as~\cite{InductionMotor-handbook}.

For permanent-magnet three-phase machines, the general structure of any physically consistent model including magnetic-saturation, saliency and other conservative effects is given by equations~\eqref{eq:PM} with magnetic Lagrangian $\LL_m$ and Hamiltonian $\HH_m$ related by~\eqref{eq:Hm}. For induction three-phase machines the physically consistent models are given by~\eqref{eq:IM} where the magnetic Lagrangian $\LL_m$ is related to the magnetic energy $\HH_m$ by~\eqref{eq:HmInd}. Such synthetic formulation of the dynamical equations is new and constitutes the first contribution of this note.  We propose here natural modifications of the standard Lagrangian to include magnetic-saturation, saliency and space-harmonics, derive the corresponding  dynamical equations and magnetic energies that could be used in the future   to construct controlled Lyapunov functions and/or storage functions  for passivity-based feedback laws.

From a control theoretic point of view we just prove here that the severe observability difficulties encountered in sensorless control and well explained in~\cite{ibarra-el-al:auto04} for the standard model resulting from the quadratic Lagrangian~\eqref{eq:IndLagrangianStandard}, remain present for models~\eqref{eq:IM}   where the magnetic Lagrangian is any  function of  the rotor angle, stator and rotor currents. Consequently,  addition of magnetic saturations, saliency and
harmonics effects, do not remove  observability issues at zero stator frequency in the sensorless case (see proposition~\ref{prop:obs}). This observability obstruction has nether been stated for models with magnetic-saturation and space-harmonics of three-phase machines and constitutes the second contribution of this notes.  Contrarily to observability, non-observability is not a generic property  and could  be destroyed by generic and small changes in the equations.   Since proposition~\ref{prop:obs} is based on the class of models derived from~\eqref{eq:PMdyn} or~\eqref{eq:IM} with arbitrary magnetic Lagrangian ${\mathcal L}_m$, we prove here that any physically consistent model of three-phase machines where the non-conservative  effects result only from voltage supply and Ohmic losses, such non-observability holds true around zero stator frequency. This means that  non-observability around zero stator frequency  is  robust to  generic and physically consistent  modifications of the equations.  As far as we know this negative and physically robust  result is new.
It  indicates that sensorless control of three-phase electrical machines around  zero stator frequency cannot be just addressed via refined physical  models but also requires   advanced and nonlinear control techniques.

In section~\ref{PM:sec} we  recall the simplest  model of   a permanent-magnet machine and its  Euler-Lagrange formulation based on the two  scalar  components of the complex stator current.  Using the  complexification procedure detailed in appendix, we show how to use  complex representation of  stator-current in  Lagrangian formulation of the dynamics. This leads us
to  the general form of physically consistent models (equation~\eqref{eq:PM}).  Finally  we obtain, just by simple  modifications of   the  magnetic  Lagrangian, physically consistent models with magnetic saturation and  saliency
effects (equation~\eqref{eq:PMgen}).  Section~\ref{Ind:sec}  deals  with induction machines and  admits the same
progression as the previous one: we start with the usual $(\alpha,\beta)$ model, describe its complex  Lagrangian
formulation,  derive  physically consistent models (equations~\eqref{eq:IM}) and  specialize them to  saturation  and
space-harmonics effects (equation~\eqref{eq:IndSatHarm} ). In section~\ref{Obs:sec}, we prove
proposition~\ref{prop:obs} that states the main observability issues of these Euler-Lagrangian models  at zero stator
frequency.   In conclusion we show how to transpose this modelling based on  complex currents associated   to a Lagrangian formulation  to complex fluxes associated to a Hamiltionian formulation with complex generalized positions and momentums.
The  appendix  details  the complexification procedure. It explains  how to derive the Euler-Lagrange equations  when some  generalized positions and velocities are treated as complex quantities. Throughout the paper, we define models in $(\alpha,\beta)$ frame, using the standard transformation from three phases frame (see, e.g.,~\cite{Leon85a}).

\section{Permanent-magnet three-phase machines }\label{PM:sec}

\subsection{The usual model and its  magnetic energy}
In the $(\alpha,\beta)$ frame (total power invariant transformation), the dynamic equations read (see, e.g.,~\cite{Leon85a}):
\begin{equation}\label{eq:PMdyn}
 \left\{
 \begin{aligned}
 &  \dotex \left(J \dot \theta\right)=
    n_p \Im\left( \left(\bar\phi e^{\jmath  n_p\theta}\right)^* \is \right) - \tau_L
   \\
 & \dotex \left( \lambda \is + \bar\phi e^{\jmath n_p\theta}  \right) = u_s-R_s \is
 \end{aligned}
 \right.
\end{equation}
where
\begin{itemize}

\item $^*$ stands for complex-conjugation, $\Im$ means imaginary part,  $\jmath=\sqrt{-1}$ and $n_p$ is the number of pairs of poles.

\item $\theta$ is  the rotor mechanical  angle,  $J$ and
$\tau_L$ are the inertia and  load torque, respectively.

\item $\is\in\CC $ is the stator current, $u_s\in\CC $ the stator voltage.

\item $\lambda = (L_d+L_q)/2$  with inductances $L_d=L_q>0$  (no saliency here).

\item The stator flux   is $\phi_s = \lambda \is+\bar \phi  e^{\jmath n_p\theta}$ with the constant  $\bar\phi >0$ representing  to the rotor flux
 due to permanent magnets.
\end{itemize}
The Lagrangian associated to this system is the sum of the  mechanical one $\LL_c$ and magnetic one $\LL_m$  defined as follows:
\begin{equation}\label{eq:energy}
    \LL_c= \frac{J}{2} \dot\theta^2, \quad \LL_m= \frac{\lambda}{2}  \left|\is+ \ib e^{\jmath n_p\theta}\right|^2
\end{equation}
where $\ib=\bar \phi / \lambda >0$ is the permanent magnetizing current.

It is well known that~\eqref{eq:PMdyn} derives from  a variational principle (see, e.g.,\cite{ortega-et-al-book98}) and thus can be written as Euler-Lagrange equations with source terms corresponding to energy exchange with the environment.
Consider the additional complex variable $q_s\in\CC$ defined by $\dotex q_s=\is$. Take  the Lagrangian $\LL=\LL_c+\LL_m$  as a  real function of the   generalized coordinates $q=(\theta,q_{s\alpha},q_{s\beta})$ and generalized velocities  $\dot q= (\dot\theta,\imath_{s\alpha},\imath_{s\beta})$:
{\small
\begin{equation}\label{eq:PMrealLagrangian}
\LL(q,\dot q) = \frac{J}{2} \dot\theta^2
   + \frac{\lambda}{2}  \left((\dot q_{s\alpha}+\ib\cos n_p\theta)^2 + (\dot q_{s\beta}+\ib\sin n_p\theta)^2\right)
\end{equation}}
Then the mechanical equation in~\eqref{eq:PMdyn} reads
$$
\dotex\left( \dv{\LL}{\dot\theta} \right) - \dv{\LL}{\theta} = -\tau_L
$$
where $-\tau_L$ corresponds to the energy exchange through the  mechanical load torque.
Similarly, the real part of complex and electrical  equation in~\eqref{eq:PMdyn} reads
$$
\dotex\left( \dv{\LL}{\dot q_{s\alpha}} \right)- \dv{\LL}{q_{s\alpha}} = u_{s\alpha} - R_s \imath_{s\alpha}
$$
and its imaginary part
$$
\dotex\left( \dv{\LL}{\dot q_{s\beta}} \right)- \dv{\LL}{q_{s\beta}} = u_{s\beta} - R_s \imath_{s\beta}
$$
since $\dv{\LL}{q_{s\alpha}}=\dv{\LL}{q_{s\beta}}=0$  and $\dot q_s= \is$. The energy exchanges here are due to the
power supply through the voltage $u_s$ and also to dissipation and irreversible phenomena due to stator resistance
represented by the Ohm law $-R_s \is$.

\subsection{Euler-Lagrange equation with complex current}
The drawback of such Lagrangian formulation is that we have to split  into real and imaginary parts the generalized complex coordinates with $q_s =q_{s\alpha} + \jmath q_{s\beta}$, ($q_{s\alpha}$ and $q_{s\beta}$ real) and velocities  $\dot q_s=\is=\dot q_{s\alpha} + \jmath \dot q_{s\beta}$, ($\dot q_{s\alpha}$ and $\dot q_{s\beta}$ real). We do not preserve  the elegant formulation of the electrical part through complex  variables and equations.

Let us apply the complexification procedure detailed in appendix  to the Lagrangian  $\LL(\theta,q_{s\alpha},q_{s\beta},\dot\theta, \dot q_{s\alpha},\dot q_{s\beta})$  defined in~\eqref{eq:PMrealLagrangian}. The complexification process only focuses on  $q_s$ and $\dot q_s =\is$ by considering $\LL$ as a function of
$(\theta,q_s,q_s^*, \dot\theta,\is,\is^*)$:
\begin{multline*}
  \LL(\theta,\dot \theta, \is,\is^*) =
   \frac{J}{2} \dot\theta^2
   + \frac{\lambda}{2}  \left(\is+ \ib e^{\jmath n_p\theta}\right) \left(\is^*+ \ib e^{-\jmath n_p\theta}\right)
   .
\end{multline*}
Following the notations in appendix, $n^c=1$ with  $q^c=q_s$, $n^r=1$ with $q^r=\theta$, $S^r=-\tau_L$ and $S^c= u_s-R_s \is$. Then according to~\eqref{eq:ELc} the usual equations~\eqref{eq:PMdyn} read
\begin{equation*}
    \dotex \left( \dv{\LL}{\dot\theta} \right) = \dv{\LL}{\theta} -\tau_L,
\quad
2 \dotex  \left( \dv{\LL}{\is^*} \right) = u_s-R_s \is
\end{equation*}
since $\dv{\LL}{q_s^*}=0$ and $ \dv{\LL}{\dot q_s^*}=\dv{\LL}{\is^*}$.

More generally, for any  magnetic Lagrangian $\LL_m$ that is a real value function of $\theta$, $\is$ and $\is^*$ and  that is $\frac{2\pi}{n_p}$ periodic versus $\theta$, we get the general  model (with saliency,
saturation, space-harmonics, ...) of three-phase permanent-magnet machine:
\begin{equation}\label{eq:PM}
   \dotex \left(J \dot \theta\right) = \dv{\LL_m}{\theta} - \tau_L,
   \quad  2 \dotex \left(\dv{\LL_m}{\is^*} \right) = u_s- R_s \is
\end{equation}
We recover the usual equation with $\phi_s = 2\dv{\LL_m}{\is^*}$ corresponding  to the stator flux and  to the conjugate momentum $p^c$ of $q^c$ as shown in appendix. According to~\eqref{eq:L2H}, the Hamiltonian $\HH$ is the sum of two energy: $\HH= \HH_c + \HH_m$. The mechanical kinetic energy $\HH_c=\frac{J}{2}\dot\theta^2$ and the magnetic energy
\begin{equation}\label{eq:Hm}
    \HH_m\left(\theta, \is, \is^*\right) =  \dv{\LL_m}{\is}\is + \dv{\LL_m}{\is^*}\is^*  - \LL_m.
\end{equation}
The   standard model~\eqref{eq:PMdyn} derives from  a  magnetic Lagrangian of the form
$
\LL_m= \frac{\lambda}{2} \left| \is  +\ib e^{\jmath n_p\theta} \right|^2
$
with $\lambda$ and $\ib$ are two positive parameters. Its corresponding magnetic energy reads
$
\HH_m =\frac{\lambda}{2} \left(|\is|^2 - \ib^2\right)$. We recover the usual  magnetic energy  $\frac{\lambda}{2} |\is|^2$ up to  the constant magnetizing energy $\frac{\lambda}{2} \ib^2$.

In~\eqref{eq:PM}, many other  formulations of $\LL_m$  are possible and  depend on
particular modeling issues. Usually, the dominant part of  $\LL_m$ will be of the form  $\frac{\bar \lambda}{2} \left|
\is + \ib e^{\jmath n_p\theta} \right|^2$ ($\bar\lambda$, $\ib$ positive constants) to which correction terms that are ``small" scalar functions of $(\theta, \is ,\is^*)$ are added.

\subsection{Saliency models}
Adding to $\LL_m$  the correction $- \frac{\mu}{2} \Re\left( \is^2 e^{-2\jmath n_p\theta} \right)$ with
$|\mu| < \lambda$ ($\Re$ means real part)  provides a simple way to represent saliency phenomena while the dominant
part of the magnetic Lagrangian (and thus of the dynamics) remains attached to $\frac{\lambda}{2} \left| \is  + \ib
e^{\jmath n_p\theta} \right|^2$. With a  magnetic Lagrangian of the form
\begin{multline} \label{eq:PMEmsaliency}
\LL_m = \frac{\lambda}{2}  \left(\is+ \ib e^{\jmath n_p\theta}\right) \left(\is^*+ \ib e^{-\jmath n_p\theta}\right)
\\
    -\frac{\mu}{4}  \left( \left(\is^*  e^{\jmath n_p\theta}\right)^2+ \left(\is e^{-\jmath n_p\theta}\right)^2\right)
\end{multline}
 where $\lambda = (L_d+L_q)/2$ and  $\mu = (L_q-L_d)/2$
(inductances $L_d>0$ and $L_q>0$), equations~\eqref{eq:PM} become ($\lambda\ib=\bar\phi$)
\begin{equation}\label{eq:PMdynsaliency}
 \left\{
 \begin{aligned}
 &  \dotex \left(J \dot \theta\right)=
    n_p \Im\left( (\lambda \is^*+ \bar\phi e^{-\jmath  n_p\theta}   - \mu \is e^{-2\jmath n_p\theta} )  \is \right) - \tau_L
   \\
 & \dotex \left( \lambda \is + \bar\phi e^{\jmath n_p\theta} - \mu \is^* e^{2\jmath n_p\theta}  \right) = u_s-R_s \is
 \end{aligned}
 \right.
\end{equation}
and we recover the usual model with saliency effect. In this case the magnetic energy  is given by:
\begin{multline}\label{eq:HmStandard}
\HH_m= \frac{\lambda}{2}  \left(|\is|^2 - \ib^2 \right)
    -\frac{\mu}{4}  \left( \left(\is^*  e^{\jmath n_p\theta}\right)^2+ \left(\is e^{-\jmath n_p\theta}\right)^2\right)
    .
\end{multline}

\subsection{Saturation and saliency models}
We can also take into account magnetic saturation effects, i.e., the fact that inductances depend on the currents. Let us assume first that only the inductances $\lambda$ and $\mu$ in~\eqref{eq:PMEmsaliency} depend on the modulus $\rho=\left|\is + \ib e^{\jmath n_p\theta}\right|$.
The  magnetic Lagrangian now reads
\begin{multline}\label{eq:PMLagrangiansaliencySaturation}
 \LL_m=
   \frac{\lambda\left( \left| \is+ \ib e^{\jmath n_p\theta}\right| \right)}{2}  \left|\is+ \ib e^{\jmath n_p\theta}\right|^2
   \\
    -\frac{\mu\left( \left| \is+ \ib e^{\jmath n_p\theta}\right| \right)}{4}  \left( \left(\is^*  e^{\jmath n_p\theta}\right)^2+ \left(\is e^{-\jmath n_p\theta}\right)^2\right)
   .
\end{multline}
 The dynamics is given by~\eqref{eq:PM} with such $\LL_m $.
Denote $\lambda^\prime = \frac{d\lambda}{d\rho}$. With
$
 \dv{\lambda}{\theta} =
  n_p\frac{\Im\left( \ib e^{-\jmath n_p\theta} \is\right)}{\left| \is+ \ib e^{\jmath n_p\theta}\right|}~\lambda^\prime$ and $
\dv{\lambda}{\is^*} =  \frac{ \is + \ib e^{\jmath n_p\theta} }{2\left| \is+ \ib e^{\jmath n_p\theta}\right|}~\lambda^\prime
$,
we get the following  dynamical model with both saliency and magnetic-saturation effects:
{\small \begin{equation}\label{eq:PMgen}
 \begin{aligned}
 &  \dotex \left(J \dot \theta\right)=
     n_p\Im\left( \left(
     \Lambda
     \left(\is^* + \ib e^{-\jmath n_p\theta}\right) - M \is e^{-2\jmath n_p\theta}\right)  \is \right) - \tau_L
   \\
 &  \dotex \left(
    \Lambda
     \left(\is + \ib e^{\jmath n_p\theta}\right) - M \is^* e^{2\jmath n_p\theta}  \right) = u_s-R_s \is
 \end{aligned}
\end{equation}
}
with $\Lambda =  \lambda+ \frac{ \left| \is+ \ib e^{\jmath n_p\theta}\right| }{2} \lambda^\prime$ and $M=\mu+ \frac{ \left| \is+ \ib e^{\jmath n_p\theta}\right| }{2} \mu^\prime$.   It is interesting to compute the magnetic energy $\HH_m$ from general formula~\eqref{eq:Hm}:
{\small \begin{multline} \label{eq:Hmsat}
\HH_m= \frac{\lambda+\left|\is+\ib e^{\jmath n_p\theta}\right|\lambda^\prime }{2} |\is|^2 -\frac{\lambda}{2}\ib^2
\\+ \frac{ \left|\is+\ib e^{\jmath n_p\theta}\right|\lambda^\prime }{4} \ib \left(\is e^{-\jmath n_p\theta}+ \is^* e^{\jmath n_p\theta} \right) \\ -
    \frac{\mu+
    \mu^\prime\frac{\is^*\left(\is+\ib e^{\jmath n_p\theta}\right)+
               \is\left(\is^*+\ib e^{-\jmath n_p\theta}\right) }
    {2\left|\is+\ib e^{\jmath n_p\theta}\right|}
    }{4}  \left( \left(\is^*  e^{\jmath n_p\theta}\right)^2+ \left(\is e^{-\jmath n_p\theta}\right)^2\right).
\end{multline}
}
Such magnetic energy formulae are not straightforward. They are not obtained by replacing $\lambda$ and $\mu$ in the standard magnetic energy~\eqref{eq:HmStandard} by $\Lambda$ and $M$ respectively.

\section{Induction three-phase machines} \label{Ind:sec}
We will now proceed as for permanent-magnet machines. Let us recall first the
  usual dynamical equations  of an induction machine  with complex stator and rotor currents. They  admit the following form
\begin{equation}\label{eq:Inddyn}
 \left\{
 \begin{aligned}
 &  \dotex \left(J \dot \theta\right)=
     n_p\Im\left( L_m \ir^* e^{-\jmath  n_p\theta}  \is \right) - \tau_L
   \\
 & \dotex \left( L_m \left(\ir+ \is e^{-\jmath n_p\theta}\right) + L_{fr} \ir \right) = -R_r \ir    \\
 & \dotex \left( L_m \left(\is+ \ir e^{\jmath n_p\theta}\right) + L_{fs} \is\right) = u_s-R_s\is
 \end{aligned}
 \right.
\end{equation}
where
\begin{itemize}

\item $n_p$ is the number of pairs of poles, $\theta$ is  the rotor mechanical  angle,  $J$ and $\tau_L$ are the inertia and  load torque, respectively.

\item $\ir\in\CC$ is the rotor current (in the rotor frame, different from  the $(d,q)$ frame) , $\is\in\CC $ the
stator current (in the stator frame, identical to the  $(\alpha,\beta)$ frame)  and $u_s\in\CC $ the stator voltage (in
the stator frame). The stator and rotor resistances are $R_s>0$ and $R_r>0$.

\item The inductances $L_m$, $L_{fr}$ and $L_{fs}$ are positive parameters with $L_{fr},L_{fs}\ll L_m$.

\item The stator (resp. rotor)  flux   is $\phi_s = L_m \left(\is+ \ir e^{\jmath n_p\theta}\right) + L_{fs} \is $ (resp. $\phi_r = L_m \left(\ir+ \is e^{-\jmath n_p\theta}\right) + L_{fr} \ir $).
\end{itemize}

\subsection{Euler-Lagrange equation with complex current}
With notations of appendix, $n^c=2$ with $q^c=(\ir,\is)$, $n^r=1$ with $q^r=\theta$, $S^c=(-R_r\ir,u_s-R_s\is)$ and $S^r=-\tau_L$. The Lagrangian associated to~\eqref{eq:Inddyn}, expressed with complex currents $\ir$ and $\is$, reads:
\begin{multline*}
  \LL(\theta,\dot \theta, \ir,\ir^*, \is,\is^*) =
   \frac{J}{2} \dot\theta^2
   \\
   +\frac{L_m}{2} \left|\is + \ir e^{\jmath n_p\theta}  \right|^2
   +\frac{L_{fr}}{2}|\ir|^2 + \frac{L_{fs}}{2}|\is|^2
   .
\end{multline*}
The first term $ \frac{J}{2} \dot\theta^2  $ represents the  mechanical  Lagrangian and the remaining sum the
magnetic Lagrangian $\LL_m$. The dynamics~\eqref{eq:Inddyn} read
\begin{multline*}
\dotex \left( \dv{\LL}{\dot\theta} \right) -\dv{\LL}{\theta} =-\tau_L,
\\
2 \dotex  \left( \dv{\LL}{\ir^*} \right) = -R_r \ir,
\quad
2 \dotex  \left( \dv{\LL}{\is^*} \right) = u_s-R_s \is,
\end{multline*}
The magnetic Lagrangian  $\LL_m$  has the following form that coincides  here  with the magnetic energy $\HH_m$:
\begin{equation}
 \label{eq:IndLagrangianStandard}
\LL_m =  \frac{L_m}{2} \left|\is + \ir e^{\jmath n_p\theta}  \right|^2
   +\frac{L_{fr}}{2}|\ir|^2 + \frac{L_{fs}}{2}|\is|^2
   .
\end{equation}
More generally physically consistent model should be obtained with a Lagrangian of the form
\begin{equation}\label{eq:IML}
    \LLIM=\frac{J}{2} \dot \theta^2 +  \LL_m\left(\theta,\ir,\ir^*,\is,\is^*\right)
\end{equation}
where $\LL_m$ is the magnetic Lagrangian expressed with the  rotor angle and currents. It is $\frac{2\pi}{n_p}$ periodic versus $\theta$.
Any physically admissible model of a three-phase induction machine reads
\begin{multline}\label{eq:IM}
   \dotex \left(J \dot \theta\right) = \dv{\LL_m}{\theta} - \tau_L,
   \\
    2\dotex \left(\dv{\LL_m}{\ir^*}\right) = -R_r\ir,
   \quad  2 \dotex \left(\dv{\LL_m}{\is^*}\right) = u_s- R_s \is,
\end{multline}
where the rotor and stator flux are given by
$$
\phi_r= 2\dv{\LL_m}{\ir^*}, \qquad \phi_s= 2\dv{\LL_m}{\is^*}
.
$$
In general the magnetic energy does not coincide with $\LL_m$. It is given by~\eqref{eq:L2H} that yields:
\begin{multline}\label{eq:HmInd}
\HH_m\left(\theta, \ir, \ir^*, \is, \is^*\right) =  \dv{\LL_m}{\ir}\ir +  \dv{\LL_m}{\ir^*} \ir^*
\\ +  \dv{\LL_m}{\is} \is +  \dv{\LL_m}{\is^*} \is^*- \LL_m
.
\end{multline}

\subsection{Saturation models}
A simple way to include saturation effects is to consider that the main  inductances $L_m$  appearing in~\eqref{eq:IndLagrangianStandard} depends on the modulus $\rho=\left| \is + \ir e^{\jmath n_p\theta}\right|$.  Thus we consider the following magnetic-saturation  Lagrangian:
\begin{multline}\label{eq:IndLagrangianSaturation}
 \LL_m=
   \frac{L_m\left(\left| \is + \ir e^{\jmath n_p\theta}\right|\right)}{2}
   \left( \is + \ir e^{\jmath n_p\theta} \right) \left( \is^* + \ir^* e^{-\jmath n_p\theta} \right)
   \\+\frac{L_{fr}}{2}\ir\ir^* + \frac{L_{fs}}{2}\is\is^*
   .
\end{multline}
Since ($L_m^\prime=\frac{dL_m}{d\rho}$)
$
 \dv{\LL_m}{\theta} = n_p
  \frac{\Im\left( \ir^* e^{-\jmath n_p\theta} \is\right)}{\left| \is+ \ir e^{\jmath n_p\theta}\right|}~L_m^\prime
$ and
$$
\dv{L_m}{\ir^*} =  \frac{ \is e^{-\jmath  n_p\theta} + \ir }{2\left| \is+ \ir e^{\jmath n_p\theta}\right|}~L_m^\prime
,\quad
\dv{L_m}{\is^*} =  \frac{ \is + \ir e^{\jmath n_p\theta} }{2\left| \is+ \ir e^{\jmath n_p\theta}\right|}~L_m^\prime
,
$$
the saturation model (formula~\eqref{eq:IM} with $\LL_m$ given by~\eqref{eq:IndLagrangianSaturation}) reads
\begin{equation}\label{eq:IndSat}
 \left\{
 \begin{aligned}
 &  \dotex \left(J \dot \theta\right)=
     n_p\Im\left( \Lambda_m \ir^* e^{-\jmath  n_p\theta}  \is \right) - \tau_L
   \\
 & \dotex \left( \Lambda_m \left(\ir+ \is e^{-\jmath n_p\theta}\right) + L_{fr} \ir \right) = -R_r \ir    \\
 & \dotex \left( \Lambda_m \left(\is+ \ir e^{\jmath n_p\theta}\right) + L_{fs} \is\right) = u_s-R_s\is
 \end{aligned}
 \right.
\end{equation}
with $\Lambda_m= L_m +\frac{\left| \is + \ir e^{\jmath n_p\theta}\right|}{2} L_m^\prime$ function of $\left| \is + \ir e^{\jmath n_p\theta}\right|$.
 We recover here  usual saturation models (see, e.g., \cite{InductionMotor-handbook}, page 428).
Notice the similarity with  permanent-magnet machines and~\eqref{eq:PMgen}. Following~\eqref{eq:HmInd}, the associated magnetic energy  reads then
\begin{multline*}
\HH_m= \frac{L_m+ \left| \is + \ir e^{\jmath n_p\theta}\right| L_m^\prime}{2}
   \left| \is + \ir e^{\jmath n_p\theta} \right|^2 \\ + \frac{L_{fr}}{2}|\ir|^2 + \frac{L_{fs}}{2}|\is|^2
   .
\end{multline*}
Notice  also that such similarity  is not complete and could be misleading:  contrarily to~\eqref{eq:IndSat} that can be derived  intuitively form~\eqref{eq:Inddyn} by replacing $L_m$ by $\Lambda_m$,  the above  magnetic energy is not provided by~\eqref{eq:IndLagrangianStandard} with $L_m$ replaced by $\Lambda_m$.

\subsection{Space-harmonics with saturation}

We can take into account space-harmonic effects by adding  their contribution to the magnetic
Lagrangian~\eqref{eq:IndLagrangianSaturation}. According to~\cite{InductionMotor-handbook}, page 298, the iron path in
general gets shorter as the harmonic order gets higher. Thus saturation effect has relatively smaller influence on the
spatial harmonics. Following~\cite{fudeh-ong:ieee83}, the Lagrangian $\LL_{\nu}$ of harmonic $\nu$ is
\begin{equation}\label{eq:harmonic}
    \LL_{\nu} = \frac{L_{\nu}}{2}\left(\is\ir^*e^{-\jmath \sigma_\nu \nu  n_p\theta}+
    \is^*\ir e^{\jmath \sigma_\nu \nu  n_p\theta}  \right)
\end{equation}
with $L_{\nu}$ a small parameter ($|L_{\nu}|\ll L_m$) and with $\sigma_\nu=\pm 1$ depending on arithmetic conditions on
$\nu$ (see~\cite{fudeh-ong:ieee83}, equations (25) to (29)).  The total
magnetic Lagrangian now reads
\begin{multline}\label{eq:IndLagrangianSaturationHarmonic}
 \LL_m= \frac{L_m\left(\left| \is + \ir e^{\jmath n_p\theta}\right|\right)}{2}
   \left( \is + \ir e^{\jmath n_p\theta} \right) \left( \is^* + \ir^* e^{-\jmath n_p\theta} \right)\\+\frac{L_{fr}}{2}\ir\ir^* + \frac{L_{fs}}{2}\is\is^*
+  \frac{L_{\nu}}{2}\left(\is\ir^*e^{-\jmath \sigma_\nu \nu  n_p\theta} +
    \is^*\ir e^{\jmath \sigma_\nu \nu  n_p\theta}  \right)
.
\end{multline}
Now the saturation model~\eqref{eq:IndSat} is changed as follows:
\small
\begin{equation}\label{eq:IndSatHarm}
 \begin{aligned}
 &  \dotex \left(J \dot \theta\right)=
     n_p\Im\left( \left(\Lambda_m e^{-\jmath  n_p\theta}
     + L_\nu\sigma_\nu \nu e^{-\jmath \sigma_\nu \nu  n_p\theta}  \right)
     \ir^*\is \right) - \tau_L
   \\
 & \dotex \left( \Lambda_m \left(\ir+ \is e^{-\jmath n_p\theta}\right) + L_{fr} \ir
 + L_\nu \is e^{-\jmath\sigma_\nu \nu  n_p\theta}
 \right) = -R_r \ir    \\
 & \dotex \left( \Lambda_m \left(\is+ \ir e^{\jmath n_p\theta}\right) + L_{fs} \is
 + L_\nu \ir e^{\jmath\sigma_\nu \nu  n_p\theta}
 \right) = u_s-R_s\is
 \end{aligned}
\end{equation}
\normalsize
Following~\eqref{eq:HmInd}, the associated magnetic energy  reads then
\begin{multline*}
\HH_m= \frac{L_m+ \left| \is + \ir e^{\jmath n_p\theta}\right| L_m^\prime}{2}
   \left| \is + \ir e^{\jmath n_p\theta} \right|^2 + \frac{L_{fr}}{2}|\ir|^2
   \\+ \frac{L_{fs}}{2}|\is|^2 + \frac{L_{\nu}}{2}\left(\is\ir^*e^{-\jmath \sigma_\nu \nu  n_p\theta}+
    \is^*\ir e^{\jmath \sigma_\nu \nu  n_p\theta}  \right)
   .
\end{multline*}
Several space-harmonics can be included in a similar way. Moreover saturation of space-harmonics can be also tackled just by choosing  $L_\nu$  as a function of $|\is +\ir e^{-\jmath  n_p\theta} |$. As far as we know such explicit models including magnetic saturation and space harmonics have never been  given.

\section{Observability issues at zero stator frequency} \label{Obs:sec}
The sensorless control case is characterized by  a load torque $\tau_L$  constant but unknown,  control inputs $u_s$ and  measured outputs $\is$. Models derived from~\eqref{eq:PM} for permanent-magnet machines (resp. from~\eqref{eq:IM} for inductions machines) can be always written in state-space form
\begin{align}
\label{eq:dyn}
\dotex X &= f(X,U), \qquad Y=h(X)
\end{align}
where $X=(\tau_L,\theta,\dot\theta,\Re(\is),\Im(\is))$ (resp.
$X=(\tau_L,\theta,\dot\theta,\Re(\ir),\Im(\ir),\Re(\is),\Im(\is))$) with $U=(\Re(u_s),\Im(u_s))$,
$Y=(\Re(\is),\Im(\is))$ and  $\dotex \tau_L=0$. A stationary regime  at  zero stator frequency corresponds then to a
steady state $(\bar X,\bar U,\bar Y)$ of~\eqref{eq:dyn} satisfying $f(\bar X, \bar U) =0$ and $\bar Y= h(\bar X)$. The
tangent linear system around this steady state is then
\begin{align}
\label{eq:dynlin}
\dotex x = A x + B u, \qquad y= C x
\end{align}
where $A=\dv{f}{X}(\bar X,\bar U)$, $B=\dv{f}{U}(\bar X,\bar U)$ and $C=\dv{h}{X}(\bar X)$.
If we assume that the linearized system~\eqref{eq:dynlin} is observable, the Kalman criteria implies  that the rank of the matrix $\left(
                                \begin{array}{c}
                                  C \\
                                  A \\
                                \end{array}
                              \right)
$ must be equal to $\dim(X)$. If it is the case,  the mapping $X \mapsto (f(X,\bar U),h(X))$ is maximum rank around $\bar X$. This maximum rank condition just means  that the set of  algebraic equations characterizing the steady-state from the knowledge of $\bar U$ and $\bar Y$,
$f(X,\bar U)=0 $ and $h(X)=\bar Y$
admits around $\bar X$ the maximum rank $\dim(X)$. Such  rank is not changed by any invertible manipulations of this set of equations characterizing the steady-state from the knowledge of  the input and output values,   $\bar U$ and  $\bar Y$. Putting the  implicit Euler-Lagrange equations~\eqref{eq:PM} and~\eqref{eq:IM} into their explicit state-space forms~\eqref{eq:dyn} involves such invertible manipulations.

For  permanent-magnet  machines described by~\eqref{eq:PM}, this set of equations  yields  the  following  mapping
$$
(\tau_L,\theta,\dot\theta,\is) \mapsto (0,\dot \theta, \dv{\LL_m}{\theta} - \tau_L,  \bar U- R_s \is , \is)
$$
where $\LL_m$ depends on $\theta$ and $\is$.  Its rank should be  maximum, i.e., equal to  $5$. This is not the case since its rank is obviously   equal to $4$.
For induction machines described by~\eqref{eq:IM}, the mapping is
$$
(\tau_L,\theta,\dot\theta,\ir,\is) \mapsto (0,\dot \theta, \dv{\LL_m}{\theta} - \tau_L, -R_r\ir,  \bar U- R_s \is , \is)
$$
where $\LL_m$ depends on  $\theta$, $\is$ and $\ir$. Its rank is  equal to  $6$  whereas the maximum rank  is $7$.   The above arguments  yield following proposition:
\begin{prop}\label{prop:obs}
Any dynamical model of permanent-magnet machines~\eqref{eq:PM} (resp. induction machines~\eqref{eq:IM}) is unobservable
around zero stator frequency regime when the measured output is the stator current $\is$ and the load torque is
constant but unknown. By unobservable we mean that:
\begin{itemize}
\item to any constant input and output $\bar u_s$ and $\bar \imath_s$ satisfying   $\bar u_s= R_s\bar \imath_s$ correspond a one dimensional family of steady states parameterized by  the scalar variable $\xi$ with
\begin{itemize}
    \item
    $\tau_L = \dv{\LL_m}{\theta}(\xi,\bar \imath_s,\bar \imath_s^*), ~ \theta =\xi,\quad  \is=\bar \imath_s $ for the permanent-magnet machines,
    \item
    $\tau_L = \dv{\LL_m}{\theta}(\xi,0,0,\bar \imath_s,\bar \imath_s^*),~ \theta =\xi, \quad \ir=0, ~\is=\bar\imath_s$ for induction machines.
\end{itemize}
\item the linear tangent systems around  such steady-states are not observable;
\end{itemize}
\end{prop}

\section{Conclusion}

The models proposed in this note, \eqref{eq:PMgen} for permanent-magnet machine and  \eqref{eq:IndSatHarm} for induction drives, are based on  variational principles and Lagrangian formulation of the dynamics. Such formulations are particularly efficient to preserve the physical insight while maintaining a synthetic view  without describing all the technological  and material  details (see~\cite{yourgrau-mandelstam-book} for an excellent and tutorial  overview of variational  principles  in physics). Extensions to network of machines and generators connected via long  lines  can also be developed  with similar  variational principles and Euler-Lagrange equations with complex currents and voltages.

In this note we have put the emphasis on currents and thus Lagrangian modelling.  Since flux variables are  conjugated to current variables, Hamiltonian modelling is also possible when fluxes are used instead of currents. For permanent-magnet machines, the Hamiltonian counterpart of Lagrangian models~\eqref{eq:PM} reads:
$$
\dotex \left(J \dot \theta\right) = -\dv{\HH_m}{\theta} - \tau_L,
   \quad  \dotex \phi_s = u_s- 2 R_s \dv{\HH_m}{\phi_s^*}
$$
where the magnetic energy  $\HH_m$  is considered as a function of  the rotor angle  $\theta$, the stator flux $\phi_s$ and its complex conjugate $\phi_s^*$. The stator current $\is$ corresponds then to $2\dv{\HH_m}{\phi_s^*}$.
For induction machines, the Hamiltonian formulation  associated to~\eqref{eq:IM} becomes
\begin{multline*}
\dotex \left(J \dot \theta\right) = -\dv{\HH_m}{\theta} - \tau_L,\\
   \dotex \phi_r = - 2 R_r \dv{\HH_m}{\phi_r^*} ,
   \quad  \dotex \phi_s = u_s- 2 R_s \dv{\HH_m}{\phi_s^*}
   .
\end{multline*}
The magnetic energy $\HH_m$ now depends on $\theta$, the rotor flux $\phi_r$ and its complex conjugate $\phi_r^*$, the stator  flux $\phi_s$ and its complex conjugate $\phi_s^*$. The rotor  (resp. stator)  current  is then given by $2\dv{\HH_m}{\phi_r^*}$ (resp. $2\dv{\HH_m}{\phi_s^*}$).
As for Lagrangian modelling, one can modify the magnetic energies of  the standard models~\eqref{eq:PMdyn} and~\eqref{eq:Inddyn} to include, for example,  magnetic-saturation or space-harmonics effects. This yields new formula expressing  $\HH_m$  as function of angle and fluxes. The corresponding flux-based models are then given by the above equations.

\appendix[Lagrangian and Hamiltonian with complex variables]
It is explained in~\cite[page 87]{cohenT-book1} how to use complex coordinates for Lagrangian and Hamiltonian systems. Here we propose a straightforward extension where the  complexification  procedure is only partial. Such extension  cannot be found directly in text-book.  In the context of electrical drives, such complexification  applies only on  electrical quantities whereas mechanical ones remain untouched.

Assume we have a Lagrangian system with generalized positions $q\in\RR^n$, $n\geq 3$ and an  analytic Lagrangian $\LL(q,\dot q)$. Let us decompose $q$ into two set of components:
\begin{itemize}
  \item the first set $q^c=(q_1,\ldots, q_{2n^c})$ with $0<2 n^c \leq n$ will be identified with $n^c$ complex numbers $q^c_k=q_{2k-1}+\jmath q_{2k}$, $k=1,\ldots, n^c$;
  \item  the second set $q^r=(q_{2n^c+1}, \ldots, q_{n})$  gathers the $n^r=n-2n^c$ components that will remain untouched and real.
\end{itemize}
Thus we can identified $q$ with $(q^c,q^r)$ where $q^c\in \CC^{n^c}$ and $q^r\in\RR^{n^r}$. Since the Lagrangian $\LL$ is a real-value and analytic  function, it can be seen as an analytic function of the complex variables $q^{c}$, $\dot q^c$, $q^{c*}$, $\dot q^{c*}$ and of the real variables $q^r$ and $\dot q^r$ ($q^{c*}$ corresponds to the complex conjugate of $q^c$). This function with be denoted by  $\LLc(q^c,q^{c*},q^r,\dot q^c, \dot q^{c*},\dot q^r)$ and is equal to $\LL(q,\dot q)$ where
{\small \begin{align*}
q&=\left( \frac{q^c_1+q^{c*}_1}{2}, \frac{q^c_1-q^{c*}_1}{2\jmath}, \ldots,
\frac{q^c_{n^c}+q^{c*}_{n^c}}{2}, \frac{q^c_{n^c}-q^{c*}_{n^c}}{2\jmath}
    , q^r_1, \ldots, q^r_{n^r} \right)
    \\
\dot q&=\left( \frac{\dot q^c_1+\dot q^{c*}_1}{2}, \frac{\dot q^c_1-\dot q^{c*}_1}{2\jmath}, \ldots,
\frac{\dot q^c_{n^c}+\dot q^{c*}_{n^c}}{2}, \frac{\dot q^c_{n^c}-\dot q^{c*}_{n^c}}{2\jmath}
    , \dot q^r_1, \ldots, \dot q^r_{n^r} \right)
    .
\end{align*}}
 Let us consider the Euler-Lagrange equations
\begin{equation}\label{eq:EL}
 \dotex \left(\dv{\LL}{\dot q_k}\right) = \dv{\LL}{q_k} + S_k, \quad k=1, \ldots, n
\end{equation}
where the $S_k$-terms  correspond to non conservative energy exchanges with the environment.  Similarly to $q$, we decompose $S=(S_1,\ldots, S_n)$ into $S^c\in\CC^{n^c}$ and $S^r\in\RR^{n^r}$: $S$ is identified with $(S^c,S^r)$.
We will  reformulate these equations with $\LLc$ and its partial derivatives. For $k=2n^c+k^\prime$ with $k^\prime=1,\ldots,n^r$, they remain unchanged $\dotex \left(\dv{ \LLc}{\dot q^r_{k^\prime}}\right) = \dv{\LLc}{q^r_{k^\prime}}+ S^r_{k^\prime}$. The two scalar equations corresponding to $k=2k^\prime-1$ and $k=2k^\prime$ with $k^\prime=1, \ldots, n^c$ yield to the following single complex equation $\dotex \left(2\dv{\LLc}{\dot q^c_{k^\prime}}\right) = 2\dv{\LLc}{q^c_{k^\prime}}+ S^c_{k^\prime}$ since we have the identities
\begin{multline*}
2 \dv{{\LLc}}{q^c_{k^\prime}}= \dv{\LL}{q_{2k^\prime-1}} -\jmath  \dv{\LL}{q_{2k^\prime}}
,\quad
2 \dv{{\LLc}}{q^{c*}_{k^\prime}}= \dv{\LL}{q_{2k^\prime-1}} +\jmath  \dv{\LL}{q_{2k^\prime}}
\end{multline*}
and similarly
\begin{multline*}
2 \dv{{\LLc}}{\dot q^c_{k^\prime}}= \dv{\LL}{\dot q_{2k^\prime-1}} -\jmath  \dv{\LL}{\dot q_{2k^\prime}}
,\quad
2 \dv{{\LLc}}{\dot q^{c*}_{k^\prime}}= \dv{\LL}{\dot q_{2k^\prime-1}} +\jmath  \dv{\LL}{\dot q_{2k^\prime}}
.
\end{multline*}
This provides the following complex formulation of the real Euler-Lagrange equations~\eqref{eq:EL}
\begin{align}
     \dotex \left(2 \dv{ \LLc}{\dot q^{c*}_k}\right) = 2\dv{\LLc}{q^{c*}_k} + S^c_k, \quad k=1, \ldots, n^c \label{eq:ELc}
     \\
      \dotex \left(\dv{\LLc}{\dot q^r_k}\right) = \dv{\LLc}{q^r_k} + S^r_k, \quad k=1, \ldots, n^r\label{eq:ELr}
      .
\end{align}
In the usual complexification procedure (\cite[page 87]{cohenT-book1}) the coefficient $2$ appearing in the above equations does not appear. This is due to our special choice $q^c_k=q_{2k-1}+\jmath q_{2k}$ instead of  the usual choice
$q^c_k=\frac{q_{2k-1}+\jmath q_{2k}}{\sqrt{2}}$. This special choice preserves  the correspondence, commonly used in electrical engineering,  between complex and real electrical quantities

Let us assume that, for each $q$, the mapping $\dot q \mapsto \dv{\LL}{\dot q}$ is a smooth bijection. Then the Hamiltonian formulation of~\eqref{eq:EL} reads
\begin{equation}\label{eq:H}
\dotex q_k=\dv{\HH}{p_k}, \quad \dotex p_k = -\dv{\HH}{q_k} + S_k, \quad k=1,\ldots,n
\end{equation}
with $\HH= \dv{\LL}{\dot q} \cdot \dot q - \LL$ and $p=\dv{\LL}{\dot q}$.   Let us decompose $p$ into  $p^c\in\CC^{n^c}$ and $p^r\in\RR^{n^r}$. Then $p^r=\dv{ \LLc}{\dot q^r}$ and $p^c = 2 \dv{ \LLc}{\dot q^{c*}}$. Simple computations yield another derivation  of the Hamiltonian from $ \LLc$ directly:
\begin{multline}\label{eq:L2H}
\HHc = \dv{ \LLc}{\dot q^c } \dot q^c  + \dv{ \LLc}{\dot q^{c*} } \dot q^{c*}
+ \dv{ \LLc}{\dot q^r } \dot q^r - \LLc
\end{multline}
where $ \HHc$ denotes the Hamiltonian $\HH$ when is a considered as a function of $(q^c,q^{c*},q^r,p^c,p^{c*},p^r)$. Then~\eqref{eq:H} becomes
{\small
\begin{align}
     &\dotex q^c_k=2\dv{\HHc}{p^{c*}_k}, \quad \dotex p^c_k = -2\dv{\HHc}{q^{c*}_k} + S^c_k,  \quad k=1,\ldots,n^c \label{eq:Hc}
     \\
     &\dotex q^r_k=\dv{\HHc}{p^r_k}, \quad \dotex p^r_k = -\dv{\HHc}{q^r_k} + S^r_k,  \quad  k=1,\ldots,n^r \label{eq:Hr}
\end{align}
}


\end{document}